\documentclass[12pt]{article}
\usepackage{amssymb, amsmath}
\usepackage{color}     

\newtheorem{thm}{Theorem}[section]

\newcommand{\io}{\int_{\{ |Q(x)| >\lambda\}} }

\renewcommand{\l}{\lambda}

\renewcommand{\o}{\omega}

\newcommand{\bb}{\begin{equation}}
\newcommand{\ee}{\end{equation}}
\newcommand{\bq}{\begin{eqnarray}}
\newcommand{\eq}{\end{eqnarray}}
\newcommand{\bqn}{\begin{eqnarray*}}
\newcommand{\eqn}{\end{eqnarray*}}

\begin{document}
\title{  Relative decay conditions on Liouville  type theorem for the steady Navier-Stokes system }
 \author{Dongho Chae\\
\ \\
Department of Mathematics\\
Chung-Ang University\\
Dongjak-gu Heukseok-ro 84\\
Seoul 06974, Republic of Korea\\
{\em email: dchae@cau.ac.kr}}
\date{}
\maketitle
\begin{abstract}
In this paper we prove Liouville type theorem for the stationary Navier-Stokes equations in $\Bbb R^3$ under  the assumptions on the relative decays of velocity,
pressure and the head pressure.
More precisely, we show that any smooth solution $(u,p)$ of the stationary Navier-Stokes equations satisfying $u(x) \to 0$ as
$|x|\to +\infty$ and the condition of  finite Dirichlet integral $\int_{\Bbb R^3} | \nabla u|^2 dx <+\infty $  is trivial, if
either  $|u|/|Q|=O(1)$ or $|p|/|Q| =O(1) $ as $|x|\to \infty$, where $|Q|=\frac12 |u|^2 +p$ is the head pressure.
  \noindent{\bf
keywords:} stationary  Navier-Stokes equations,  Liouville type theorem \end{abstract}

\section{Introduction}
 \setcounter{equation}{0}

We are concerned on the Liouville type problem for the  stationary Navier-Stokes equations in $\Bbb R^3$.
\bb\label{sns}
\left\{ \aligned u\cdot \nabla u &=-\nabla p +\Delta u,\\
\nabla \cdot u&=0
\endaligned
\right.
\ee
Here, $u=(u_1 (x), u_2 (x),   u_3 (x))$ is a vector field in $\Bbb R^3$, and $p=p(x)$ is a scalar field.  
We equip \eqref{sns}  with the uniform decay condition at spatial infinity, $|u(x)|\to 0$ as $|x|\to +\infty,$
which means more precisely 
\bb\label{bc1}
\lim_{R\to+\infty} \sup_{|x|>R} |u(x)|=0.
\ee
Obviously  $(u,p)$ with $u=0$ and $p=$constant  is a trivial solution to \eqref{sns}-\eqref{bc1}. An important question is if there is other nontrivial solution.
This uniqueness problem, or  equivalently the  Liouville type problem is currently a hot subject of study in the community of mathematical fluid mechanics.
In general we impose an additional  condition, the finiteness of the Dirichlet integral,
\bb\label{dir}
\int_{\Bbb R^3} |\nabla u|^2 dx <+\infty.
\ee
A solution $u$ of the stationary Navier-Stokes equations, satisfying the finite Dirichlet integral condition \eqref{dir}, is called D-solution, or Leray solution, following 
the work due to Leray\cite{ler}.  The Liouville type problem for D-solution of \eqref{sns}-\eqref{bc1}  is written explicitly  in \cite[Remark 9.6, p.147]{gal}. After the statement of the problem 
it is shown that   under the assumption $u\in L^{\frac{9}{2}} (\Bbb R^3)$  for  the solution the solution  we have $u=0$. The two dimensional Liouville type
problem is solved in \cite{gil}, while in the $n$-dimensional case with $n\ge 4$, the problem is easily solved  by Galdi \cite{gal}.  Thus, the Liouville type problem for D-solution has been wide open only in the three dimensional case.
There are numerous partial results(see e.g. \cite{cha1,cha11, cha3, cha2, cham, koch, koro, koz, lyu, ser1, ser2, ser3} and the references therein),  deducing the triviality of solution  to \eqref{sns}-\eqref{dir} under various sufficient conditions. 
Among others the author of this paper obtained in \cite{cha1} that under the condition together with \eqref{bc1}
$
\int_{\Bbb R^3} |\Delta u|^{\frac65} dx <+\infty
$
the solution of \eqref{sns}  is trivial, $u=0$ on $\Bbb R^3$. This is the first condition with the same scaling as \eqref{dir}.
Another sufficient condition  with the same scaling as \eqref{dir} deduced in \cite{cha11} is
$$
\int_{\Bbb R^3}|\nabla \sqrt{|Q|}|^2 dx= \frac14\int_{\Bbb R^3}\frac{|\nabla Q|^2}{|Q|}  dx<+\infty,
$$
where $Q=\frac12 |u|^2 +p$ is the head pressure.\\
\ \\
In this paper we prove that the solution of \eqref{sns} is trivial  under the conditions \eqref{bc1} and \eqref{dir} under the extra assumptions on the relative decay assumptions on the 
velocity, pressure and the head pressure.  More precisely we shall prove the following.

 \begin{thm} Let $(u,p)$ be a smooth solution of \eqref{sns} satifsying the conditions \eqref{bc1} and \eqref{dir}, and let 
 $Q=\frac12 |u|^2 +p$ is its head pressure.
Suppose that  either
\bb\label{co1}
\lim_{R\to +\infty }  \sup_{\{ |x|>R\}}   \frac{|u|}{|Q|^\frac12}  <+\infty,
\ee
or 
\bb\label{co2}
\lim_{R \to +\infty }  \sup_{\{ |x|>R\}}    \frac{|p| }{|Q|}  <+\infty.
\ee
Then, $u=0$ on $\Bbb R^3$.
\end{thm}
{\bf Remark 1.1 } Observe that the ratios in \eqref{co1} and \eqref{co2}  are scaling invariant. We also note that a  lower estimate part of the ratio in \eqref{co2}, namely
\bb\label{rem}
\lim_{R\to +\infty}  \sup_{\{ |x|>R\}}   \frac{|p| }{|Q|} \ge 1
\ee
holds true.   Indeed,   from the fact that $Q\le 0$ on $\Bbb R^3$(see \eqref{hyp} below), and $p=-\frac12 |u|^2 +Q \le 0$  we find
 $\frac{|p| }{|Q|} = \frac{1}{|Q|} ( |Q|+ \frac12 |u|^2) \ge  1$ on $\Bbb R^3,$ which provides us with \eqref{rem}.


\section{Proof of the Theorem}
 \setcounter{equation}{0}
\noindent{\bf Proof of Theorem 1.1  } 
Given smooth solution $(u,p)$ of \eqref{sns},  we  set the head pressure $Q = p +\frac12 |u |^2$ and the vorticity  $\o$ is identified as collection of independent element of the skew-symmetric part of the matrix $\nabla u$, namely $\frac12 \nabla u- \frac12 (\nabla u)^\top  $. Then,  it is well-known that the following holds
\bb\label{ns}
\Delta Q -u \cdot \nabla  Q =|\o|^2.
\ee
Indeed,  multiplying the first equation of \eqref{sns} by $u$, we have
\bb\label{nns1}
u\cdot \nabla Q=\Delta \frac{|u|^2}{2} -|\nabla u|^2.
\ee
Taking divergence of \eqref{sns}, we are led to
\bb\label{nns2}
0=\Delta p+ \rm{Tr}(\nabla u (\nabla u)^\top) .
\ee
Adding \eqref{nns2} to \eqref{nns1}, and observing $ |\nabla u|^2-   \rm{Tr}(\nabla u (\nabla u)^\top)=|\o|^2$,  we obtain \eqref{ns}.
We also recall that it is known (see e.g. Theorem 6.1, p.105, [2]) that  the condition \eqref{bc1} together with \eqref{dir} implies  that
 there exists a constant $p_0$ such that 
 \bb\label{th1}
 p(x) \to p_0  \quad \text{as}\quad |x|\to +\infty.
 \ee
Therefore,  redefining  $Q-p_0$ as a new head pressure,  denoted by the same $Q$, we may assume  that
\bb\label{bc2}
Q(x)\to 0 \qquad \text{as}\quad |x|\to+\infty.
\ee

  We shall now prove $Q=0$ on $\Bbb R^3$ by a  contradiction argument. Let us suppose on the contrary $Q\not=0$ on $\Bbb R^3$.  
 In view of  \eqref{bc2},  applying the maximum principle  to  \eqref{ns}, we observe $Q\leq 0$ on $\Bbb  R^3$.
Moreover, by the maximum principle again, either $Q(x) = 0$  for all $x\in \Bbb R^3$, or $Q(x) < 0$ for all 
$x\in \Bbb R^3$.  Indeed, any point  $x_0 \in \Bbb R^3$ such that $Q(x_0) =0$ is a point of local maximum, which is not allowed unless $Q \equiv 0$ by the maximum principle.  
Therefore, the assumption   $Q\not=0$ on $\Bbb R^3$ implies that  
\bb\label{hyp}
Q(x)<0 \quad \forall x\in \Bbb R^3.
\ee
Let us set 
$ m=\sup_{x\in \Bbb R^3} |Q(x)| >0.$
We first claim that  
\bb\label{21}
\int_{\{ |Q(x)| >\lambda \}} u\cdot \nabla  Q f(|Q|) dx =0
\ee
for any continuous function $f$ on $\Bbb R$, and for almost every $\lambda \in [0, m]$. 
Indeed,
\begin{align}
&\int_{\{ |Q(x)| >\lambda \}} u\cdot \nabla  Q f(|Q|) dx= -\int_{\{ |Q(x)| >\lambda \}} u\cdot \nabla \left( \int_0 ^{|Q(x)|} f(s) ds \right) dx\cr
&\qquad \qquad = -\int_{\{ |Q(x)| =\lambda \} } u\cdot \nu \left( \int_0 ^{|Q(x)|} f(s) ds \right)   dS \cr
&\qquad\qquad =  -\left( \int_0 ^{\lambda} f(s) ds \right)  \int_{\{ |Q(x)| =\lambda \}} u\cdot \nu \,dS\cr
&\qquad  \qquad = -\left( \int_0 ^{\lambda} f(s) ds \right)  \int_{\{ |Q(x)| >\lambda \}} \nabla \cdot  u \,\,dx=0,
  \end{align}
  where $\nu$ denotes the outward unit normal vector on the boundary of the domain 
  $$\Omega_\lambda:=\{x\in \Bbb R^3\, |\,  |Q(x)| >\lambda \},
  $$
   and we used Sard's theorem to  justify the use the divergence theorem on the domain  $\Omega_\lambda $,   whose boundary $\partial \Omega_\lambda$, is a smooth level surface for each $\lambda \in (0, m)$ expect possibly measure zero set of points.
  Integrating \eqref{21} over  $\Omega_\l$, and using \eqref{21} for $f(s)\equiv1$, we are led to
  \bb\label{22}
  \io \Delta Q=\int_{ \{ |Q(x)|=\l \}} |\nabla Q(x) |dS =  \io |\o|^2 dx.
  \ee
  Given  $\theta\in (0,1) $ we introduce a smooth, non-decreasing  
  cut-off function $\eta =\eta_\l ; [0, \infty) \to [0, 1] $ by
 \bb
 \eta (s) = \left\{ \aligned & 0 &\quad \text{if} \quad  &0\le s< 1-\theta, \\
     & 1 &\quad \text{if} \quad  &s\ge 1, 
  \endaligned \right.
  \ee
  together with  
$
0\le \eta^\prime (s) \le  \frac{2}{\theta }
$
for all $s\in (0, \infty)$.
Then, we set 
\bb
\varphi _\l (x)= \eta \left( \log \left( \frac{|Q(x)|}{\l} \right) \right) .
\ee
Note that $\varphi \in C^\infty _c ( \Bbb R^3)$ for almost every  $\l \in [0,\frac{m}{e^2} )$, and
\bb
\varphi_\l (x)= \left\{ \aligned & 0 &\quad \text{if} \quad  &0< |Q(x)| <  \l e^{1-\theta} , \\
 & 1 &\quad \text{if} \quad  & |Q(x)| \ge   \l e,
  \endaligned \right.
  \ee
and therefore, $\varphi_\l (x) \to 1$ as $\l \to 0$ for almost everywhere in $\Bbb R^3$.
We claim that for $\theta \in (0, 1)$
\bb\label{23}
\int_{\{ \theta \l <|Q(x)|  <  \l \}} \frac{|\nabla Q|^2}{|Q|} dx \le \left(  \log \frac{1}{\theta}\right) \int_{\{ |Q(x)| >\theta \l\}}  |\o|^2 dx.
\ee
Indeed, using the co-area formula and \eqref{22},  we estimate
\begin{align}
&\int_{\{\theta \l <|Q(x)|  <  \l \}} \frac{|\nabla Q|^2}{|Q|} dx = \int_{\theta \l} ^\l  \frac{1}{s}\left( \int_{ \{ |Q(x)| =s\}} |\nabla Q(x)| dS \right) ds\cr
&\qquad = \int_{\theta \l} ^\l  \frac{1}{s}\left( \int_{ \{ |Q(x)| >s\}}|\o|^2 dx \right) ds\le \int_{\theta \l} ^\l  \frac{ds}{s} \int_{ \{ |Q(x)| >\theta\l \}}|\o|^2 dx \cr
&\qquad =\left(  \log \frac{1}{\theta}\right)\int_{\{ |Q(x)| >\theta \l\}}  |\o|^2 dx,
\end{align}
as claimed.   We observe that this implies
\bb\label{23a}
\int_{\{e^{ 1-\theta}  \l <|Q(x)|  <  e^\l \}} \frac{|\nabla Q|^2}{|Q|} dx \le  \theta  \int_{\{ |Q(x)| >e^{ 1-\theta}  \l \}}  |\o|^2 dx.
\ee

Multiplying  \eqref{sns} by $u$, we find
  \bb\label{25}
  \Delta |u|^2 -2 (u\cdot \nabla )Q = 2 |\nabla u|^2.
  \ee
  Next, multiplying \eqref{25} by $\varphi _\l (x)$,  integrating it over $\Bbb R^3$,  and integrating by part, we obtain
  \begin{align}\label{26}
  \int_{\Bbb R^3}  |\nabla u|^2 \varphi_\l  dx&= - \frac12 \int_{\Bbb R^3}  \frac{ \nabla |u|^2 \cdot \nabla Q }{|Q| } \eta^\prime \left( \log \left( \frac{|Q(x)|}{\l} \right) \right)  dx \cr
   &\qquad -\int_{\{ |Q(x)| > e^{1-\theta} \l \}}   (u\cdot \nabla )Q   \, \varphi_\l   dx \cr
   & =I_1+I_2.
  \end{align}
    We assume now \eqref{co1}.   From  \eqref{21}  it follows immediately that $I_2=0$. Using \eqref{23a}, we  estimate $I_1$ as follows
  \begin{align}\label{27}
  I_1&\le \frac{2}{\theta} \int_{\{e^{(1-\theta)}  \l <|Q(x)|  < e \l \}}  \frac{|u|}{|Q|^\frac12} \frac{|\nabla Q|}{|Q| ^{\frac12}}      |\nabla u| dx\cr
  &\le\frac{2}{\theta}\sup_{\{e^{(1-\theta)}  \l <|Q(x)|  < e \l \}}  \frac{|u|}{|Q|^\frac12}   \left( \int_{\{e^{(1-\theta)}  \l <|Q(x)|  < e \l \}} \frac{ |\nabla Q|^2}{|Q|  } dx  \right) ^{\frac12}  \times\cr
  &\qquad\qquad\times 
   \left( \int_{\{e^{(1-\theta)}  \l <|Q(x)|  < e \l \}} \  |\nabla u|^2 dx \right) ^{\frac12}\cr     
  &\le\frac{2}{\sqrt{\theta} }   \sup_{\{e^{(1-\theta)}  \l <|Q(x)|  < e \l \}}   \frac{|u|}{|Q|^\frac12} \left(\int_{\{ |Q(x)| >e^{(1-\theta)}  \l \}} |\o|^2 dx \right)^{\frac12}
  \times \cr
   &\qquad\qquad \times \left( \int_{\{e^{(1-\theta)}  \l <|Q(x)|  < e \l \}} \  |\nabla u|^2 dx \right) ^{\frac12}.   \cr \end{align}
   We fix $\theta \in (0,1)$. Let $B_R=  \{ x\in \Bbb R^3 \, |\, |x|<R\} $.   
   From \eqref{bc2} we find that for  any $R>0$ there exists $\l=\l_R$ such that 
   $$ B_R  \subset  \{ x\in \Bbb R^3 \, |\, |Q(x)| > e \l\}  \subset    \{ x\in \Bbb R^3 \, |\, |Q(x)| > e^{1-\theta}  \l\},    $$
   and
   $$ 
   \l_R \to 0\quad\text{as}\quad R\to +\infty.
   $$
  From this observation we have from \eqref{26} and \eqref{27}, replacing $\l$ by $\l_R$, that
   \begin{align}\label{28}
   \int_{B_R}  |\nabla u|^2  dx     &\le\frac{2}{\sqrt{\theta} }   \sup_{\{    |x| >R \}} 
     \frac{|u|}{|Q|^\frac12} \left(\int_{\Bbb R^3}  |\o|^2 dx \right)^{\frac12}
   \left( \int_{\{     |x|>R \}} \  |\nabla u|^2 dx \right) ^{\frac12}.  
   \end{align}
Passing $R\to \infty $ on the both sides of \eqref{28},     and applying  the Lebesgue dominated convergence theorem,   
   we obtain $\int_{\Bbb R^3} |\nabla u|^2 dx =0$, and $u=$constant on $\Bbb R^3$. Combing this with the boundary condition \eqref{bc1},  we  find $u=0$
    on $\Bbb R^3$. Regarding  the condition \eqref{co2},  we  just observe 
    $$ \left( \frac{|u|}{|Q|^{\frac12} }\right) ^2 \le   \frac{ 2|Q| +2|p|}{|Q|} =     2+ 2 \frac{|p|}{ |Q|}, $$
    and
    $$
    \frac{ |p|}{|Q|}\le \frac{|Q|+  \frac{|u|^2}{2} }{|Q|} \le 1 +\frac12 \left( \frac{|u|}{ |Q|} \right)^2.
    $$
Hence, the conditions \eqref{co1} and  \eqref{co2}  are equivalent to each other.  $\square$\\
  \ \\

        $$\mbox{\bf Acknowledgements}$$
This research was  supported partially by NRF grants 2016R1A2B3011647. 
 
\end{document}